\title{ On the separation profile of infinite graphs }

\author{Itai Benjamini \and Oded Schramm \footnote{Oded died while solo climbing Guye Peak in Washington State on September 1, 2008.}\and \'Ad\'am Tim\'ar}

\documentclass[12pt,naturalnames]{article}
\usepackage{amsmath}
\usepackage{amsthm}
\usepackage{amsfonts}
\usepackage{graphicx}
\newif\ifhyper\IfFileExists{hyperref.sty}{\hypertrue}{\hyperfalse}
\ifhyper\usepackage{hyperref}\fi

\newif\ifdraft
\drafttrue
\numberwithin{equation}{section}
\numberwithin{figure}{section}

\newtheorem{theorem}{Theorem}
\numberwithin{theorem}{section}
\newtheorem{corollary}[theorem]{Corollary}
\newtheorem{lemma}[theorem]{Lemma}

\newtheorem{proposition}[theorem]{Proposition}
\newtheorem{conjecture}[theorem]{Conjecture}

\newtheorem{question}[theorem]{Question}

\newcommand{\R}{\mathbb{R}}

\newcommand{\Z}{\mathbb{Z}}
\newcommand{\N}{\mathbb{N}}

\def\H{\mathbb{H}}

\def \E {{\bf E}}

\def \proof {{ \medbreak \noindent {\bf Proof.} }}

\def \_reg {\rightarrow_{\bf reg}}
\def \sep{cut}
\def\maxdeg/{\Delta}

\def\bintree{T_2}

\begin{document}
\maketitle

\bigskip

\begin {abstract}
Initial steps in the study of inner expansion properties
of infinite Cayley graphs and other infinite graphs, such as hyperbolic ones, are taken, in a flavor similar to the well-known Lipton-Tarjan
$\sqrt{n}$ separation result for planar graphs.
Connections to
relaxed versions of quasi-isometries are explored, such as regular and
semiregular maps.
\end {abstract}

\section{Introduction}


The best-known separation type result is the Lipton-Tarjan \cite{LT} theorem,
claiming that there is a way of removing $O(\sqrt{n})$ vertices (and the edges reaching them)
from an $n$-vertex planar graph, so that each of the remaining connected components has at most $n/2$ vertices
(see also \cite{AST1} \cite{AST2} \cite{MTT1} \cite{MTT2} for other proofs
and generalizations).

Rather than considering finite graphs, we will look at infinite graphs,
and consider the separation function. The
{\it separation function} at $n$ is the supremum over all subgraphs of size $n$,
of the number of vertices needed to be removed from the subgraph, in order to cut it
to connected pieces of size at most $n/2$. We are interested in separation functions up to constant factors. See the exact definitions
in Section~\ref{s.def} below.

The separation function (or profile) is a natural coarse geometric invariant of
infinite graphs and path metric spaces.
In this work we would like to view the separation function
as such an invariant, in the family of  invariants like isoperimetric profiles
and  volume growth.

One question we are interested in is,
which separation functions are possible for transitive graphs or Cayley graphs?
What is the separation function of a given Cayley graph? E. g., what is it for groups of
intermediate growth?  Our results below will give examples of Cayley graphs with
finite separation (free groups), logarithmic separation (Planar hyperbolic groups),
$n^{(d-1)/d}$-separation, $d \in \N$ ($\Z^d$, $\H ^{d+1}$),  $n/ \log n$-separation (Product of free groups).

We don't know of any example of a vertex transitive graph with separation $n^{\alpha}, 0 < \alpha < 1/2$.
It is maybe of interest to look at separation function for natural families of groups,
e.g. groups generated by finite automata (\cite{N}), lamplighter groups, random groups, Kazhdan groups... The latter could be natural candidates for groups with linear separation, if there is any.

\begin{question}
\label{linear}
Is there linear separation, or $n/ \log n$
is the largest growth possible?
\end{question}
In which groups do balls or the minimizers of the isoperimetric problem or
convex sets admit the largest separation among all sets
of a fixed size?


Understanding separation is useful in figuring out
the partial order given by {\it regular maps} (see Subsection~\ref{s.def} below) on the collection of
all spaces, or rather, on the most familiar ones.

In \cite{BScp} it was shown that bounded degree transient planar graphs admit non-constant bounded Dirichlet harmonic functions.
Gil Kalai asked  in 1994 whether square root separation can replace planarity? The answer is negative, see \cite{BScp}.
By a {\it comb} we mean the $\Z^2$ grid when all the edges parallel to $x$-axis but not on the $x$-axis removed.
A $ \mbox{{\it comb}} \times \Z$ is an example of a transient graph with square root separation that
does not admit non-constant bounded harmonic functions (in particular, non-constant harmonic Dirichlet functions), and thus is not planar \cite{BScp}.

\begin{question}
\label{sp}
Does spectral radius $< 1$  plus an additional separation condition imply the
existence of non-constant bounded harmonic functions?
\end{question}

Indeed, finite separation transient graphs admit non-constant bounded harmonic functions: for any such bounded size cut separating the graph to two components if only one component were transient we would get a transient simple random walk escaping to infinity via infinitely many cut sets of bounded size, which is impossible, yet having two transient components connected via a finite set implies the existence of non-constant bounded harmonic functions. What about a weaker separation condition?

In \cite{BK}  it is shown that if to
each level of a binary tree the edges of an expander sequence are added, the result has no non-constant bounded harmonic functions. This explains that some separation condition for 
Question~\ref{sp} is necessary.


In Section~\ref{s.FS} we prove that finite separation is the same as bounded treewidth (see definition there), and we present a theorem about
the structure of infinite graphs with finite separation.

Then we study the separation of products of graphs, Section~\ref{s.Pro}.
First, a bound on the separation
function of a product is given. Later it is shown that
for a regular tree $T$,
$$
sep_{T \times T}(n) \asymp  {n \over \log n},
$$
where $\asymp $ means up to constant factor. For functions $f$ and $g$, we write $f=O(g)$ or $g=\Omega (f)$ to denote that there exists a $c>0$ such that $f(x)\leq c g(x)$ for every $x$.
In Section~\ref{s.hyp}, on Gromov-hyperbolic graphs (or hyperbolic graphs for short), we open with establishing the separation function for $\H^d$ followed by a gap theorem, showing that
for hyperbolic graphs the separation function is either a constant, or is growing
at least logarithmically.

A regular map between two graphs is a map that increases distances at most by some linear function, and such that there is a uniform bound on the number of preimages of a vertex (see Subsection~\ref{s.def} for the precise definition).
Given a regular map, the separators for a constant neighborhood of the image can be pulled back, thus the separation function
is monotone non decreasing under regular maps (see Lemma~\ref{newlemma}).
Section~\ref{s.reg} discusses regular maps.

In the last section
it is proved that there is no regular map
from $\Z^2$ to the $\Z_2$ lamplighter over $\Z$.
The notion of a semi-regular map is introduced and discussed. Asymptotic dimension is monotone under semi-regular maps, hence this seems to be the right type of function for its study (the way that regular maps are suitable for separation).
\medskip

Many open problems are scattered along the paper.


\subsection{Some definitions} \label{s.def}

\noindent
{\bf Definition.} Let $G$ be a graph.
Suppose $S \subset G$ is given. Let {\bf $\sep (S)$} denote the infimum of the  sizes of
subsets $C_S$ of $S$, so that the largest connected component of $S - C_S $ has size
smaller than $|S|/2$. For such a $C_S$ we will say (with a slight sloppiness in terminology), that $C_S$ {\it separates} $S$.

The {\bf separation function} $sep_G(x):\R \to \R$ is
$$
sep_G(x)= \sup_{S \subset G, |S|=x} \sep (S).
$$
\medskip

\noindent
{\bf Remark:} The separation function can be defined in a wider context of path metric spaces,
Riemannian manifolds, in particular.
For simplicity of the exposition we stick to graphs. When the separation of a manifold is discussed,
the reader can replace the manifold by a rough isometric graph, defined as follows.
\medskip

\noindent
{\bf Definition.} Let $(X,d_{X}),(Y,d_{Y})$ be metric spaces, and let
$\kappa < \infty$. A $\kappa$-rough isometry (or $\kappa$ quasi-isometry) $f$ from $X$ to $Y$ is a
(not necessarily continuous) map $f:X \rightarrow Y$ such that
$$
\kappa^{-1} d_X(x_1,x_2)- \kappa \leq d_Y(f(x_1),f(x_2)) \leq \kappa d_X(x_1,x_2)+\kappa
$$
\noindent
holds for all $x_1,x_2 \in X$, and for every $y_1 \in Y$ there is some
$x_1 \in X$ such that
$$
d_Y(y_1,f(x_1)) \leq \kappa.
$$
\noindent
If such an $f$ exists, we say that $X$ and $Y$ are {\bf roughly isometric}.
\medskip

It is straightforward to check that roughly isometric graphs have the same separation. In particular, the separation function is a group-invariant for finitely generated groups (it does not depend on which Cayley graph is chosen). However, there is a coarser equivalence relation that still preserves the separation function. This equivalence relation can be defined by the existence of {\it regular maps} between the two graphs in both directions (and graphs will turn out to have monotone increasing separation under regular maps).

\medskip
\noindent
{\bf Definition.} Let $(X,d_{X}),(Y,d_{Y})$ be metric spaces, and let
$\kappa < \infty$. A (not necessarily continuous) map $f:X \rightarrow Y$
is  $\kappa$-regular if the following two conditions are satisfied.

\noindent
$(1)$  $d_Y(f(x_0),f(x_1)) \leq \kappa(1+d_X(x_0,x_1))$ holds for every
$x_0,x_1 \in X$, and
$(2)$  for every open ball $B=B(y_0,1)$ with radius $1$ in $Y$, the inverse
image $f^{-1}(B)$ can be covered by $\kappa$ open balls of radius $1$ in $X$.

\noindent
A {\bf regular} map is a map which is  $\kappa$-regular for some
finite $\kappa$. Write $X \rightarrow_{\bf reg} Y$ if there is a regular map from
$X$ to $Y$.

It is easy to check that if there is a rough isometry between bounded degree graphs $X$ and $Y$, then there is a regular map from $X$ to $Y$ (and also to the other direction, by the symmetry of being roughly isometric -- which is not apparent from the definition). Hence being roughly isometric implies the existence of regular maps from $X$ to $Y$ and from $Y$ to $X$.
However, the existence of a regular map from $X$ to $Y$, and one from $Y$ to $X$, does not imply that $X$ and $Y$ would be rough isometric. For example, consider two copies of $\Z^2$ ``glued" along the  axis $\{(0,n)\}$ for $X$, and $\Z^2$ for $Y$; see \cite{BScp}.

Regular maps were studied in a somewhat different context by David and Semmes \cite{DS}.
One use of the separation function is as an obstruction for the existence of a regular map
from one graph to another. This is so because if there is a regular map from $X$ into $Y$, then
$ \sep_X(x) \leq C sep_Y(x)$ for some $C < \infty$.

Let $\sep^c(G)$ denote the minimal number of vertices that are
necessary to separate the finite graph G into $c$ times
smaller pieces.  Define
$sep^c_G(x):= \sup_{S \subset G, |S|\leq x} \sep^c (S)$. (In particular, $sep^{1/2}_G(x)=sep_G(x)$.) Then
\begin{equation}
\label{newlemmaeq}
sep^c_G(x)\asymp sep_G(x).
\end{equation}

To see this, suppose $c<1/2$, and consider
some graph of size $x$. After we cut it into pieces of size at most half of the original one, we can repeat the procedure for the smaller pieces, and iterating if necessary, we can get to graphs of size at most $c$ times that of the original one. Now, in the constantly many steps of the procedure, each separating set used along the way is bounded by  $sep^{1/2}_G(x)$, hence their total size is at most a constant times $sep^{1/2}_G(x)$.  (When $c>1/2$, do the same thing, but starting from a $c$-separating set.)

\begin{lemma}
\label{newlemma}
Let $X$ and $Y$ be graphs with a uniform upper bound $d$ on their degrees, and suppose that there is a $\kappa$-regular map $f$ from $X$ to $Y$. Then $ sep_X(x) = O( sep_Y(x))$.
\end{lemma}

\noindent
{\bf Proof:}
Let $A\subset X$ be an arbitrary set of vertices that induces a connected graph. Define $A'\subset Y$ as the $2\kappa$-neighborhood of $f(A)$. By $\kappa$-regularity, $A'$ is connected, and it has size at most $|A|d^{2\kappa}$. Let $S'$ be a minimal subset of $A'$ that separates it into pieces $C'_1,\ldots, C'_m$, each of size at most $A'/(2\kappa d^{2\kappa})$.

Let $S_0$ be the $2\kappa$-neighborhood of $S'$ in $Y$. Then $S:=f^{-1} (S_0)$ has size $\leq \kappa d^{2\kappa} |S'|=O(\sep_Y (|A'|))=O(\sep_Y (|A|))$, using (~\ref{newlemmaeq}). 
Denote by $C_i$ the preimage of $C'_i$ by $f$.

We claim that $S$ is a separating set in $A$ between $C_1\setminus S, C_2\setminus S,\ldots , C_m\setminus S$. Suppose not: then there is a path $P$ in $A\setminus S$ between some $C_i\setminus S$ and $C_j\setminus S$ ($i\not =j$).  Then $f(P)$ is disjoint from $f(S)=S_0$, thus
the $2\kappa$-neighborhood $P'$ of $f(P)$ in $A'$ does not intersect $S'$. Since $P'$ is connected, this shows that some vertex of $f(C_i\setminus S)\subset C_i'$ and some vertex of $f(C_j\setminus S_0)\subset C_j'$ is connected by a path inside $A'\setminus S'$. This contradicts the fact that $C_i$ and $C_j$ are different components of $A'\setminus S'$.

We have seen that $|S|=O(\sep_Y (|A|))$, and just seen that $S$ is a separating set in $A$ between $C_1\setminus S, C_2\setminus S,\ldots , C_m\setminus S$. On the other hand, $|C_i\setminus S|\leq |C_i|\leq \kappa |C_i'|\leq |A|/2$.
\qed

We know and used three ways to rule out existence of regular maps between spaces:
separation, Dirichlet harmonic functions, and growth. See \cite{BScp} (where $\kappa$-regular maps are  called $\kappa$-quasimonomorphisms). A fourth way is to use that asymptotic dimension is monotone under regular (and more generally, semi-regular) maps; see Section~\ref{s-reg}.

Separation and growth are monotone with respect to regular maps.
Flows with finite energy can be pushed back and forth with a regular map.
One should look  for additional invariants.
A partial motivation to study separation is to try to figure out
the partial order given by regular maps on the collection of
all spaces, or rather, on the familiar ones.

\section{Finite separation and regular maps}\label{s.FS}

A graph admits the finite separation property $(FS)$ iff $sep_{G}(n)$ is a
bounded function.
Trees, and graphs that are roughly isometric to trees, are examples
of graphs with finite separation.
Yet these are not the only graphs with finite
separation. Consider an infinite Sierpi\'nski graph, which we define as some reasonable limit (say, local convergence, with root chosen to be always one of the three extremal vertices)
 of the sequence of finite Sierpi\'nski graphs.
This is an example of a graph
with $sep_{G}(n) \leq 3$ which is not roughly isometric to a tree, as it
contains arbitrary large cycles. In this
section we will try to understand the structure of graphs with
finite separation.
\medskip

\begin{theorem}
\label{twsep}
 If a bounded degree $G$ has finite separation then $G$ admits a regular map
to the 3-regular tree.
\end{theorem}

{\bf Remark:}
It is easy to check that every locally finite graph has a regular map to a tree if the tree can have arbitrarily large degrees (because the uniform bound on the number of preimages in the definition of a regular map can be ignored by ``blowing up" vertices to large enough stars).

The proof for Theorem~\ref{twsep} will proceed by showing that finite separation implies bounded treewidth.

\def\tw{{\rm tw}}
\medskip
\noindent
{\bf Definition}
Let $G$ be a finite graph, $T$ be a tree, and consider a family ${\cal V}=(V_t)_{t\in T}$ such that $V_t\subset V(G)$ for every $t$. We say that $(T, {\cal V})$ is a {\it tree-decomposition} of $G$ if the following hold
\begin{enumerate}
\item $\cup_{t\in V(T)} V_t =V(G)$;
\item for every $e\in E(G)$ there is a $t$ such that both endpoints of $e$ are in $V_t$;
\item for every $x\in V(G)$, the set $\{t\in V(T)\, :\, x\in V_t\}$ induces a connected subgraph of $T$.
\end{enumerate}
The {\it width} of the tree decomposition is $\max_{t\in T} |V_t|-1$.

The {\it treewidth} of $G$, denoted by $\tw (G)$, is the minimum of the width of all tree decompositions of $G$.

An important property of tree decompositions is that for every edge $e=\{x,y\}$ of $T$, the set $V_x\cap V_y$ is a separating set between $\cup_{t\in C_1} V_t\setminus (V_x\cap V_y)$ and $\cup_{t\in C_2} V_t\setminus (V_x\cap V_y)$ , where $C_1$ and $C_2$ are the two components of $T\setminus e$.

The following theorem was proved by Robertson and Seymour \cite{RS}. See Theorem 12.4.4 in \cite{D} for a proof, and for more details on treewidth.

\noindent
\begin{theorem}
\label{grid-minor}
For every $m$ there is a $k$ such that a graph of treewidth at least $k$ contains an $m$ by $m$ square grid as a minor.
\end{theorem}

\def\tw{{\rm tw}}

\noindent
\begin{lemma}
\label{separationtw}
For any finite graph $G$, $\sep (G)\leq \tw (G)+1$. Conversely, for every $c$ there exists a $k$ such that if $\tw (G)\geq k$ then $\sep (G) \geq c$.
\end{lemma}

\noindent
{\bf Proof:}
To see the first claim, let $(T, {\cal V})$ be a tree decomposition of $G$ of width $\tw (G)$. Choose $x\in V(T)$ so that the maximum of $|\cup_{t\in C} V_t\setminus V_x|$ over all components $C$ of $T\setminus x$ is minimal. It is easy to check that this maximum is at most $|V(G)|/2$; on the other hand, this is the maximum size of components in $G\setminus V_x$. Hence $V_x$ separates, and the first part is proved.

For the other assertion, let $c$ be any positive integer.
Choose $m\geq (4c^4 +1)^{1/2}(c^2+1)$ to be an integer, and $k$ as in Theorem~\ref{grid-minor} for this $m$. We will show that if $G$ has treewidth at least $k$, then it cannot be separated by $c$ vertices.

By assumption and Theorem~\ref{grid-minor}, such a $G$ contains an $m$ by $m$ grid $M$ as a minor. Let $G'$ be a minimal subgraph of $G$ that can be contracted to a graph isomorphic to $M$. For each vertex $x$ of $M$, let $W_x$ be the set of vertices of $G'$ that got contracted into $x$. Suppose that the vertices $x_1,x_2,\ldots , x_{m^2}$ of $M$ are ordered so that
$$|W_{x_1}| \geq |W_{x_2}|\geq \ldots \geq |W_{x_{m^2}}|.
\label{e1}
$$
Consider the set $X:=\{x_1,x_2,\ldots , x_{4c^4}\}$ in $M$. Subdivide $M$ into (at least $4c^4+1$
many)
 $c^2+1$ times $c^2+1$ pairwise disjoint subgrids (that is, subgraphs that are grids). At least one of these subgrids does not intersect $X$; let $M_0$ be a such subgrid. Let the vertices on the boundary of $M_0$ be $y_1,y_2,\ldots , y_{4c^2}$. Our plan is to choose pairwise disjoint paths $P_1,\ldots , P_{2c^2}$ from $X$ to $\{y_1,\ldots ,y_{4c^2}\}$ in $M$. Using these, for each $y_i$ that is the endpoint of some $P_j$, we will replace $W_{y_i}$ by a $\tilde W_{y_i}$ that is still connected, and has ``many" vertices (and the same number for each such $y_i$). Then we will show that $(\cup \tilde W_{y_i}) \bigcup (\cup_{y\in M_0} W_y)$
is a subgraph of $G$ with no separation of size at most $c$. We will explain this construction in more detail later.

It is easy to check, using the isoperimetry of the square grid, that there is no separating set of size $\leq 2c^2$ between $X$ and $\{y_1,\ldots , y_{4c^2}  \}$ in $M_0$. Hence, by the the max-flow-min-cut theorem, $P_1,\ldots , P_{2c^2}$ can be chosen. There is no restriction in assuming that their endpoints in $\{y_1,\ldots , y_{4c^2}\}$ are respectively $y_1,\ldots, y_{2c^2}$.
For every $i\in\{1,\ldots , 2c^2\}$, define $G_i$ to be the subgraph of $G'$ induced by $\cup_{v\in P_i} V(W_v)$. The $G_i$ are connected (by the fact that $G_i$ can be contracted to $P_i$), and pairwise disjoint. For each $i\in\{1,\ldots , 2c^2\}$, let $W'_{y_i}$ be a connected subgraph of $G_i$ that contains $W_{y_i}$ and has size $\mu := \min \{|W_{x_1}|,\ldots , |W_{x_{4c^4}}|  \}$. The choice of such $W_{y_i}'$s  is possible since $|W_{y_i}|\leq\mu$ and  $|W_{x_i}|\geq \mu$
(note that $|W_{y_i}|\leq |W_{x_j}|$ for every $i, j\in \{1,\ldots , 4c^4\}$, by the choice of $M_0$), and by the fact that the endpoint of each $P_i$ is some element of $\{x_1,\ldots, x_{4c^4}\}$.
For every $y\in M_0\setminus \{y_1,\ldots, y_{2c^2}\}$, let $W_y':=W_y$.
Finally, let $H$ be the subgraph of $G'$ induced by $\cup_{y\in M_0} W_{y}'$. We claim that $H$ is not separable by $\leq c$ vertices.

To see this, consider some separating set $S$ of size at most $c$ in $H$. Define $S_{M_0}$ as $\{y\in M_0\, :\, S\cap W_y'\not =\emptyset\}$. Then $|S_{M_0}|\leq c$, and hence $M_0\setminus S_{M_0}$ has a component $C$ of size at least $|M_0|-c^2=(c^2+1)^2-c^2$. In particular, at least $c^2$ of the $y_1,\ldots , y_{2c^2}$ are in $C$. The subgraph of $H$ induced by $\cup_{y\in C} W_y'$ is clearly connected, and by the last sentence, it has at least $c^2\mu$ vertices. On the other hand,
$\cup_{y\in M_0\setminus C} W_y'\leq \mu |M_0\setminus C|\leq c^2\mu$. This shows that the complement of $\cup_{y\in C} W_y'$ in $H\setminus S$ has size at least $c^2\mu$, while the union of all other components and $S$ has size at most $c^2\mu$. Thus $S$ does not cut $H$ to pieces of proportion $<1/2$. Since $S$ was arbitrary, this finishes the proof.
\qed



\bigskip

\def\stw{{\rm stw}}

\noindent
{\bf Proof of Theorem~\ref{twsep}:}

\begin{lemma}\label{bddtree} Let $T$ be a bounded degree tree.
Then there is a regular map $T \_reg \bintree$.
\end{lemma}

We leave the proof of the lemma as an (easy) exercise  to the reader.
\medskip

The {\it strong treewidth} of a finite graph $H$, denoted $\stw (H)$, is defined as follows. Take a so-called {\it strong tree decomposition}, which is defined to be a pair $(T,{\cal V})$, ${\cal V}=(V_i)_{i\in T}$, where $T$ is a tree, and $(V_i)_{i\in T}$ is a {\it partition} of $V(H)$, and for every edge $\{x,y\}\in E(H)$, $x\in V_i$, $y\in V_j$, then either $i=j$ or $\{i,j\}\in E(T)$. The strong treewidth, $\stw (H)$, is defined as the minimum of $\max _{i\in T} |V_i|$ over all strong tree decompositions.

It is shown in \cite{BE} that for any family of bounded degree graphs, having bounded treewidth is equivalent to having bounded strong treewidth (see Theorem 95 in \cite{Bo}). Take now the family $\Gamma$ of all finite connected subgraphs of $G$. From Lemma~\ref{separationtw} we know that finite separation of $G$ implies that there exists a $c$ such that $\tw (H)\leq c$ for every $H\in \Gamma$. Hence there exists a $c'$ such that $\stw (H)\leq c'$ for every $H\in\Gamma$. Given an $H$, take a strong tree decomposition $(T,{\cal V})$. (We may assume that every $V_i\in{\cal V}$ is nonempty.) Then the map $\phi : H\rightarrow T$, $x\mapsto i$ where $x\in V_i$ is 1-regular. Observe that the maximal degree of $T$ is at most $c'd$, otherwise there were adjacent vertices $s$, $t$ in $T$ such that there is no edge between $\cup_{i\in C_1} V_i$ and $\cup_{i\in C_2} V_i$, where $C_1$ and $C_2$ are the components one gets after deleting the edge $\{s,t\}$ from $T$. This would contradict the connectedness of $H$ (since every $V_i$ is nonempty). So, for every $H\in\Gamma$, there is a $\phi$ 1-regular map into a tree of maximal degree at most $c'd$. A standard compactness argument gives then that there is a 1-regular map from $G$ to a $c'd$-regular tree. Since the composition of regular maps is regular, Lemma~\ref{bddtree} finishes the proof.
\qed

\section{Products}\label{s.Pro}

How does separation behave
under products? Below is  a lower bound that might be
also the right upper bound.  At least it is tight for the
products $\Z^{d_1+d_2} = \Z^{d_1} \times \Z^{d_2}$.

\begin{theorem}
\label{prod}
For $c=(7/8)^{1/2}$ and any finite graphs $G,H$
$$
\sep(G \times H) \asymp  \min(|H| \sep^c(G), |G| \sep^c(H)),
$$
where $\asymp$ means up to constant factors which are independent
of $G$ and $H$.
\end{theorem}

\noindent
{\bf Proof:} Take a finite graph $G$, and let $K$ be a subset of $V(G)$.
How can we tell if $K$ separates $G$?  Here is how.
Pick two vertices $v,u$ in $G$ at random independently and uniformly.
Let $p$ be the
probability that they are not in the same component of
$G \setminus K$ (including the case when at least one of them is in $K$).  Then $K$ separates only if $p > 1/2$.

Now consider some set of vertices $K$ in $G \times H$.
Let $(g_1,h_1)$ and $(g_2,h_2)$ be two randomly independently,
uniformly chosen points in $G \times H$.  Let
$p_1$ be the probability that $(g_1,h_1)$ and $(g_1,h_2)$ are
in different components of ${g_1} \times H-K$, let $p_2$ be the
probability that $(g_1,h_2)$ and $(g_2, h_2)$ are in different
components of $G \times {h_2}-K$, and let $p$ be
the probability that $(g_1,h_1)$ and $(g_2,h_2)$ are in different
components of $G \times H-K$.  Clearly, $p \leq p_1 + p_2$.
If $K$ separates, then $p > 1/2$, so at least one of $p_1,p_2$
is greater than $1/4$.  Suppose that it is $p_1$.
That implies that for at least $(1/8) \times |G|$ choices of $g_1$ in $G$
the probability that $h_1$ and $h_2$ are in different
components of ${g_1} \times H-K \cap ({g_1} \times H)$ is at least $1/8$.
So for those $g_1$ we have $K \cap ({g_1} \times H)$ separating ${g_1} \times H$
(for  $c=(7/8)^{1/2}$, that is, the components after deletion have sizes at most $c$ times the original size).
Hence, for such $g_1$ the
cardinality of $K \cap ({g_1} \times H$) is at least $\sep ^c(H)$.
Hence the cardinality of $K$ is at least $C|G| \sep ^c(H)$ with some $C$.
This proves the direction
$$
   \sep(G \times H) > C \min(|H| \sep^c(G), |G| \sep^c(H)).
$$
The other direction is obvious.
\qed

Recall (\ref{newlemmaeq}).  Using that, Theorem~\ref{prod} provides us with a lower bound for the separation function of product graphs.
However, it does not settle the question regarding the separation
of a product of infinite graphs, because
it does not provide us with an upper bound
for the separation of a finite graph $F \subset G \times H$ that is not
a product.

\noindent
\begin{theorem}
\label{newtheorem} Let $G_1$ and $G_2$ be two (possibly infinite) graphs. Then
$$
sep_{G \times H} (N) = \Omega(  \max_k\min({N\over k} sep_G (k), k sep_H (N/k))).
$$
In particular,
$$sep_{G\times\Z} (m^2/ sep_G (m) )=\Omega ( m).
$$
\end{theorem}

\noindent
{\bf Proof:}
The first half is straightforward from Theorem~\ref{prod} (using also (\ref{newlemmaeq})). For the second assertion, let $N=m^2 / sep_G (m)$. The numbers whose minimum we take on the right side are equal (up to constant factor) if $k=m/ sep_G (m)$, hence the maximum on the right, given by this case, is $m$.
\qed

\noindent
\begin{corollary}
\label{calcul}
$$sep_{\Z^d}(n)\asymp n^{(d-1)/d},$$
$$sep_{\H^2\times\H^2}(n)=\Omega (n^{1/2}/\log n),$$
$$sep_{\R\times \H^2}=\Omega (n\log n)^{1/2}.$$
\end{corollary}

\noindent
{\bf Proof:}
The lower bound for $sep_{\Z^d} (n)$ follows from the finite separation of $\Z$ using induction and Theorem~\ref{newtheorem}. The upper bound is given in
Proposition~\ref{hypspace} (see the details there), and by a similar argument one can get the upper bound for $\R \times \H^2$ (also using Theorem~\ref{TxT}).

For $\H^2\times \H^2$, we
use that the separation of $\H^2$ is $\log n$ by Theorem~\ref{newtheorem}. By symmetry we have
$ \max_k\min({n\over k} sep_{\H^2} (k), k sep_{\H^2} (n/k))=n^{1/2} \log n^{1/2}\asymp n^{1/2} \log n$.

Finally,
the lower bound for $sep_{\R\times\H^2}$ is straightforward from the $G\times \Z$ part of Theorem~\ref{newtheorem}.
\qed

The previous lower bound on $sep_{\H^2\times\H^2}$ will be improved in Theorem~\ref{hrtt}.

\medskip


\begin{conjecture}
$sep_{G \times \Z}(n)$ is always equal to the lower bound
in Theorem~\ref{newtheorem}.

This means that up to constants
$sep_{G \times \Z}(m^2 / sep_G(m)) \asymp m$.
\end{conjecture}

The value of $sep_{\H^2 \times \R}(n)$ is a particular case of that,
the lower bound was given in Corollary~\ref{calcul}.

\medskip

We are bugged by not being able to settle the question of the separation
of a graph of the form $G \times \Z$, where the separation function of $G$ is known.
If the worst case is a product subgraph of $G \times \Z$ then
$sep_{G \times \Z}(n^2/sep_G(n)) \asymp n$
would be the answer.
%
We can prove this for the case when $G$ has finite separation. To obtain $sep_{G \times \Z}(n) \asymp \sqrt n$, take a regular map from $G$ to a regular tree $T$ (using Theorem~\ref{twsep}). This defines a regular map from ${\Z \times G}$ to $\Z\times T$. The separation of the latter is $\sqrt{n}$; see Lemma 7.2 in \cite{BScp}, hence we have the $\sqrt n$ upper bound. For the lower bound, apply Theorem~\ref{prod}.

Denote by $T$ the binary tree, we have:

\noindent
\begin{theorem}
\label{TxT}
$$
sep_{T \times T}(n) \asymp  {n \over \log n}.
$$
\end{theorem}

\noindent
{\bf Remark:} This implies an upper bound of $n/ \log n$ for the product of any
two bounded valence trees (using Theorem~\ref{bddtree}), and a lower bound for $\H^2 \times \H^2$.

The proof is similar to that of Lemma 7.2 in \cite{BScp}.
\medskip

\noindent
{\bf Proof:} Lower bound first. A useful representation of $T \times T$ is
as follows. Consider
sequences 
that are finite, consist of $0$'s and $1$'s, and have a decimal point.
Two of these are neighbors if one is obtained from the other
by adding a digit either on the right or on the left.
  We will use this model for $T \times T$.
Let $B_k$ be the set of points in $T \times T$ that are at
distance $k$ from the root (which is ".").  These are just the
length $\leq k$ binary sequences with a dot at some place inside.  Suppose
that $W$ is a separating set for $B_k$.  Pick two length $2k$ binary sequences at
random, with the dot in the center of each.  Let these be called
$a$ and $b$.  At random, pick a $k$-length subword of $a$, and denote
it $a'$, and pick a $k$-length subword of $b$, denoted $b'$.  Then $a'$
and $b'$ are just independent uniformly selected random vertices in $B_k \setminus B_{k-1}$.
Let $a_1$ be the left half of $a$, that is, the part of $a$ to the
left of the dot.  This is a point in $B_k$.  Let $a_2$ be the
right half of $a$.  Same for $b=b_1.b_2$, all have length $k$.  Consider the following path
from $a'$ to $b'$ through $a_2$ and $b_1$.  Delete the leftmost digit in $a'$.
Then add to the right the next digit of $a_2$.  Do so repeatedly,
until you reach $a_2$.  Then delete the rightmost digit, and add
on the left the digit of $b_1$.  Continue that until you reach $b_1$.
Then delete digits on the left and add on the right the digits
of $b_2$, until you reach $b'$.  In this path you see length
$k$ and $k-1$ subwords of the words $a$, $b_1.a_2$, and $b$.  With
probability at least $1/4$, say, you must pass through the separating
set $W$ in this process.  Since $a$, $b_1.a_2$ and $b$ are equally
distributed, the probability that $b$ contains a subword in $W$
(of length $k$ or $k-1$) is at least $1/12$.  There are $4^k$ different
possible $b$.  Each word in $W$ with length $k$ or $k-1$ is contained in
$2^k$ or $2 \times2^k$ different $b$'s.
Hence the size of $W$ is at least $2^k/24$.
But $|B_k| \asymp k \times 2^k$.  So $|W| > C |B_k|/ \log(|B_k|)$.

Now, the upper bound.  Let $S$ be a finite set of points in $T \times T$.
We shall find a separating set for $S$ of size $C n/ \log n$, where
$n=|S|$.  Let $p_1$ and $p_2$ be the projections into the factors.  Let
$D(p_1(v)), D(p_2(v))$ be the distance from the root in the respective factors.  Let $m_1$
be the median of $D(p_1(v))$ on $S$, and let $m_2$ be the median
of $D(p_2(v))$ on $S$.
Let $A_1(m)$ be the set of vertices in $S$ such that $D(p_1(v)) =m_1$,
and similarly for $A_2(m)$.  If the size of $A_1(m_1)$ is
less than $100 n/ \log n$, then we can take $A_1(m_1)$ as the separating
set.  Otherwise, let $k^+_1$ be the least $k \geq m_1$ such
that $|A_1(k^+_1)| < 100 n/ \log n$, let $k^-_1$ be the greatest
$k \leq m_1$ such that $|A_1(k^-_1)| < 100 n/ \log n$.  Similarly,
define $k^-_2$ and $k^+_2$.  The separating set is
$A = A_1(k^-_1) \cup A_1(k^+_1) \cup A_2(k^-_2) \cup A_2(k^+_2)$.
The size of $A$ is obviously less than $400 n / \log n$.
We just need to check that $A$ separates well.
Let $K$ be a component of $S \setminus A$, and let $v$ be some vertex in $K$.
If $D(p_1(v))< k^-_1$, then every vertex $w$ in $K$ satisfies
this inequality.  Because $m_1 \geq k^-_1$ is the median for
$D(p_1(z))$, this clearly implies that $|K| \leq n/2$.
So we may assume that $D(p_1(v)) \geq k^-_1$ holds
for every vertex $v$ in $K$.  Similarly, one gets
$D(p_1(v)) \leq k^+_1$,  $k^-_2 \leq D(p_2(v)) \leq k^+_2$.
A connected component of ${t \in T : k^- \leq D(t) \leq k^+}$
has cardinality exactly $2^{k^+ - k^- -1}$.  Hence it follows that
$ |K| \leq 2^{k^+_1 - k^-_1} \times 2^{k^+_2 - k^-_2}$.
But note that $k^+_j - k^-_j \leq \log(n)/100 +1$ for $j=1,2$,
because for each $k$ in the range $k^-_j < k < k^+_j$
we have $|A_j(k)| \geq 100 n/ \log n$.
This gives
    $|K| \leq 4 \times 2^{\log n/50} < n/2$.
Hence $A$ is a separating set.
\qed

The proof  shows also that $sep_{T \times T \times T}(n) \asymp n/ \log n$,
and the same for any $T \times T \times T... \times T$.

\medskip

\section{Hyperbolic graphs} \label{s.hyp}

Before proceeding to the study of more general hyperbolic graphs, we determine the separation function of the hyperbolic space $\H ^d$. This follows from the method in \cite{MTT1}, \cite{MTT2}, which can also be used to give a proof for the separation of $\R^d$.

\begin{proposition}
\label{hypspace}
For $d=2$, $sep_{\H ^d}\asymp \log n$, and for $d\geq 3$, $sep_{\H ^d}\asymp n^{(d-2)/(d-1)}$.
For $d\geq 1$, $sep_{\Z ^d}\asymp n^{(d-1)/d}$.
\end{proposition}

\proof
Suppose $d\geq 3$. The lower bound follows from the fact that $\R ^{d-1}$ embeds isometrically into $\H ^d$.

For the upper bound, we adapt the proof of \cite{MTT1} to our setting.

Take a connected fundamental domain $Q$ of $\H^d$ by some group of isometries acting
co-compactly and properly discontinuously on $\H^d$. The translates of $Q$ by this group give a
tiling ${\cal T}$ of $\H^d$; denote the corresponding dual graph by $G$. For a vertex $x\in G$, let $\tau (x)$ be the corresponding tile.
Clearly $G$ is transitive and roughly isometric to $\H^d$. We want to
show that $G$ has separation $n^{(d-2)/(d-1)}$ if $d\geq 3$, and $\log
n$ if $d=2$.

So let $H$ be a finite subgraph of $G$, and ${\cal T}|_H$ be the union of the
corresponding set of tiles in $\H^d$. Denote by $o$ the center of mass
of ${\cal T}|_H$; note that any hyperplane $P$ through $o$ cuts ${\cal T}|_H$ into
two pieces of equal volume. In particular, $\{x\in H: \tau (x)$ intersects $P\}=:S$ separates $H$. So we want
to bound $|S|$.

Choose $P$ randomly and uniformly of all hyperspaces through $o$. Look
at ${\cal T}$ in the corresponding Poincar\'e ball model; we may assume that
$o$ is the origin (and then $S$ is the intersection of a Euclidean
hyperspace with the open unit ball). Let the tile containing the origin be $t$. If $d\geq 3$, the expected number of tiles in ${\cal T}|_H$ intersected by $P$ can be estimated as follows. The set of all tiles in ${\cal T}$ at graph-distance $k$ from $t$ is of order $c^{(d-1)k}$ with some $c>1$; the set of those intersecting $P$ is of  order $c^{(d-2)k}$ (here we are using the fact that the tiles all have the same volume, hence the same $c$ works for both cases). Hence, if $m_k$ is the number of tiles in ${\cal T}|_H$ intersecting $P$ and at graph distance $k$ from $t$, we obtain $\E (|S|)=O(\sum_{k=1}^\infty m_k c^{(d-2)k}/c^{(d-1)k})=O(\sum m_k c^{-k})$. This is maximized when $H$ is a ball, in which case we have $\E(|S|)=O(c^{(d-1)\log n}c^{-\log n})$
(where the base of the logarithm is $c^{d-1}$). Hence, there exists a $P$ such that the corresponding $S$ has $O(c^{(d-2)\log n})=O(n^{(d-2)/(d-1)})$ elements - and this is what we wanted to prove.

The proofs for $\H^2$ and for $\Z^d$ proceed similarly, so we omit the details.

\qed

From this proof technique of \cite{MTT1},
one also gets that a graph which can be sphere packed in $\R^d$ such that the spheres in the packing have bounded ratios, have separation at most $O( n^{(d-1)/d})$. This implies e.g. that
$\Z^{d+1}$ cannot be sphere packed in $\R^d$. For more on this, see \cite{BSpack}.
\medskip

\noindent
{\bf Remark:} In \cite{KL} the following is proved. Say that a graph $G$ has growth rate $d$ if every ball of radius $r$ in $G$ contains at most $r^d$ vertices. Then there is an injective graph homomorphism from $G$ to
$\Z^{d\log d}_\infty$, where $\Z^{d\log d}_\infty$ is the graph on $\Z^{d\log d}$ where two vertices are adjacent if each of their coordinates differ by at most one. Together with the method of the previous proof, we can conclude that $sep_G (n)=O (n^{(d\log d -1)/(d\log d)})$ for such a $G$.

\bigskip

\noindent
In the standard hyperbolic plane $\H^2$, triangles do not have the same properties
as in the Euclidean plane. For instance, in the Euclidean plane, in any large
isosceles triangle, the midpoint of the hypotenuse is far away from the
other two sides. This cannot happen in hyperbolic space. That observation
led E. Rips to the following definition.
\medskip

\noindent
{\bf Definition.} Let $G=(V,E)$ be a graph. Given three vertices $u,v,w \in V$,
pick geodesics between any two to get a geodesic triangle. Denote the
geodesics by $[u,v],[v,w],[w,u]$. Say the triangle is $\delta$-thin if for any
$v' \in [u,v]$

$$
\min(d(v',[w,u]),d(v',[v,w])) \leq \delta
$$

\medskip
\noindent
$G$ is said to be $\delta$-hyperbolic if there is some $\delta \geq 0$, such
that all geodesic triangles in $G$ are $\delta$-thin.

Hyperbolic groups were introduced by Gromov ~\cite{G2}.
They are among the central objects in geometric group theory.

We do not have any general upper bounds on the separation of hyperbolic
graphs.
Still, we have the following gap theorem.


\begin{theorem}
\label{hypqit}
Let $G=(V,E)$ be a graph which is $\Delta$-hyperbolic
and has maximal degree $M<\infty$.
Let $N$ be some integer.
There is a $c>0$, which depends only on $M,N$ and $\Delta$,
such that if
\begin{equation}
\label{sublog}
sep_G(n) < c\log n\qquad \text{for all }n>N,
\end{equation}
then $G$ is roughly-isometric to a tree.
\end{theorem}

\begin{lemma}
\label{lemdiam}
Let $G=(V,E)$ be a graph which is $\Delta$-hyperbolic
and has maximal degree $M<\infty$, and let
$v_0$ be some vertex of $G$.

Let us write $u\sim v$ whenever $u,v$ are two
vertices of $G$ with $d(u,v_0)=d(v,v_0)$ and
there is a sequence of vertices $u=u_1,u_2,\dots,u_k=v$ in
$G$ with $d(u_j,v_0)=d(u,v_0)$ and
$d(u_j,u_{j+1})\leq2\Delta$ for each $j=1,\dots,k-1$.

Let $N$ be some integer.
There are $D,c>0$, which depend only on $M,N$ and $\Delta$,
such that if (\ref{sublog}) holds, then
$d(u,v)<D$ whenever $u\sim v$.
\end{lemma}

\proof
For each vertex $v\neq v_0$, let $\sigma_1v$ be one of the
neighbors of $v$ that is closer to $v_0$, and set $\sigma_1 v_0=v_0$.
Inductively, define $\sigma_{n+1}v=\sigma_1\sigma_nv$,
and set $\sigma(v)=\{v,\sigma_1v,\sigma_2v,\dots,v_0\}$.
Note that this is a geodesic from $v$ to $v_0$,
and $\sigma(u)\subset \sigma(v)$ holds whenever $u\in \sigma(v)$.

Now assume that $u_1,u_2,\dots,u_k$ is a sequence of vertices such that
$d(u_j,v_0)=d(u_1,v_0)$ and $d(u_j,u_{j+1})\leq 2\Delta$
for each $j=1,\dots,k-1$.
Also assume that (\ref{sublog}) holds for some $c$ satisfying
\begin{equation}
\label{cchoice}
c^{-1} > (2\Delta+3) \log M,
\end{equation}
and for some $N$.
We need to show that $d(u_1,u_k)<D$ for some $D=D(\Delta,M,N)$.

It follows from (\ref{cchoice}) that
\begin{equation}
\label{need}
c \log 11 + c \log Q + c Q \log M < \frac{Q-\Delta}{2\Delta + 3}.
\end{equation}
holds for all sufficiently large $Q$.  So let $Q$ be an integer large enough
to satisfy (\ref{need}) and
$$
10 M^Q>\max\{N,2\Delta\}.
$$
Also set
$$
L= 10 M^Q,\qquad R=d(u_1,v_0).
$$

It is clear that we may assume $R>L+Q$,
for otherwise $D= 2 (L+Q)$ suffices.
Let $t$ be the largest integer in $\{1,\dots,k\}$
such that the cardinality of
$A=\{\sigma_{L+Q} u_j:j=1,\dots,t\}$ is at most $Q$.
See figure below.
Our plan is to prove that $|A|<Q$.  The definition of
$t$ then shows that $t=k$.
As we shall see, the diameter of $A$ is bounded by $2\Delta |A|$.
Hence it follows from $t=k$ and $|A|<Q$ that
$d(u_1,u_k)\leq D$ with $D=2\Delta Q + 2 L + 2 Q$
(consider the path from $u_1$ to $u_k$ constructed by following
$\sigma(u_1)$ to $A$, taking the shortest path
to $A\cap \sigma(u_k)$, and then following
$\sigma(u_k)$ to $u_k$), and the lemma follows.
So the main task is to show that $|A|<Q$.
\medskip


For every $w\in A$ let $x_w$ be one of the $u_j$'s in
$\{u_1,\dots,u_t\}$ so that $w = \sigma_{L+Q} x_w$,
and let $s_w$ be the segment of $\sigma(x_w)$ between
$x_w$ and $w$;
that is, $s_w=\{\sigma_jx_w:j=1,2,\dots,L+Q\}$.

Let $U$ be the set of vertices of $G$ that are at distance
at most $Q$ from $A$ or that are in one of the segments
$s_w, w\in A$, and let $H$ be the graph obtained by
restricting $G$ to $U$.
By construction, the segments $s_w, w\in A$ are disjoint, and
consequently,
$$
|U| \geq |A| (L+Q) > |A| L = 10 |A| M^Q.
$$
On the other hand, a ball of radius $Q$ in $G$ has
clearly no more than $M^Q$ vertices, and it follows that
$$
|U| \leq |A| ( L + M^Q) = 11 |A| M^Q.
$$
{}From (\ref{sublog}) and $10 M^Q>N$
it follows that there is a set of vertices $U_0\subset U$
with
\begin{equation}
\label{uzbd}
|U_0| \leq c \log |U| \leq c \log ( 11 |A| M^Q ),
\end{equation}
such that every component of $H-U_0$ has less than
$6 |A| M^Q$ vertices.

We are trying to prove that $|A| < Q$.
By construction $|A|\leq Q$, so assume that $|A|=Q$.
Therefore, from (\ref{uzbd}) and (\ref{need}) it follows that
\begin{equation}
\label{qq}
|U_0| \leq c \log 11 + c \log Q + c Q \log M <
 \frac{Q-\Delta}{2\Delta + 3}\leq Q/3.
\end{equation}

Let $A'$ be the set of $w\in A$ such that the segment
$s_w$ does not meet $U_0$.
Because these segments are disjoint, and $|A|=Q$,
it follows from (\ref{qq}) that
$$
|A'|>2Q/3.
$$

Let $r$ be an integer in the range $[R-L-Q,R-L-\Delta]$,
and let
$$
V_r = \{v\in V: |d(v,v_0)-r|\leq\Delta\}.
$$
We claim that there is a path whose vertices
are in $V_r\cap U$ that intersects each $s_w, w\in A$.
For each $j=1,\dots,t-1$, let $y_j$ be the vertex on $\sigma(u_j)$
that has distance $r$ from $v_0$.  Consider the triangle
in $G$ formed by taking $\sigma(u_j),\sigma(u_{j+1})$ and
a shortest curve $\alpha_j$ joining $u_j$ and $u_{j+1}$.
Since $G$ is $\Delta$-hyperbolic, there is a vertex
$z_{j+1}$ in $\alpha_j\cup \sigma(u_{j+1})$ whose distance
to $y_j$ is at most $\Delta$.
Because the length of $\alpha_j$ is at most $2\Delta$,
the distance from $\alpha_j$ to $v_0$ is at least
$R-\Delta$.  On the other hand,
$$
d(z_{j+1},v_0)\leq d(z_{j+1},y_j)+d(y_j,v_0)\leq \Delta+r < R-\Delta,
$$
and it follows that $z_{j+1}\notin\alpha_j$.
Hence $z_{j+1}\in \sigma(u_{j+1})$, and there is a path
$\beta_j$ in $U$ from $y_j$ to $z_{j+1}$ with all vertices in $V_r$ (since $d (z_{j+1}, v_0)<R-L$).
Let $\gamma_{j+1}$ be the arc of $\sigma(u_{j+1})$ connecting
$z_{j+1}$ with $y_{j+1}$.
It is now clear that the the union of all the arcs $\beta_j$ and
$\gamma_j$ is a curve whose vertices are in $V_r\cap U$
that intersects each $s_w,w\in A$.
This construction with $r=R-Q-L$ also shows that
the diameter of $A$ is bounded by $2 \Delta |A|$, as we have promised.

Let $m$ be the largest integer smaller than $(Q-\Delta)/(2\Delta+2)$.
Set $r_j=R-Q-L+2 j (\Delta+1)$, for $j=0,1,\dots,m$,
and note that the sets $V_{r_j}$ are disjoint.
Since
$$
|U_0|<m+1
$$
follows from (\ref{qq}),
there must be some such $r_j$ so that $U_0$ does not meet $V_{r_j}$.
{}From the above it then follows that there is a connected
component of $H-U_0$ that intersects each $s_w,w\in A$.
Consequently, there is a connected component $X$ of $H-U_0$ that
contains each $s_w,w\in A'$.  The number of vertices of $X$
is at least
$$
|A'| L \geq 2 Q L /3 = 20 |A| M^Q/3 > 6|A| M^Q.
$$
This contradicts the definition of $U_0$, and
the contradiction establishes $|A|<Q$ and the lemma.
\qed

\begin{lemma}
\label{tree}
Let $G$ be a $\Delta$-hyperbolic graph, let
$v_0$ be some vertex in $G$, and let $\sim$ be the equivalence
relation defined in Lemma~\ref{lemdiam}.
Let $V_1$ be the set of vertices of $G$ whose distance from
$v_0$ is divisible by $3\Delta+1$, and let
$\tilde V=V_1/\sim$ be the set of equivalence classes of $\sim$
in $V_1$.
Let $\tilde G=(\tilde V,\tilde E)$ be defined by
letting $[\tilde v,\tilde u]\in\tilde E$  whenever
for some $v\in \tilde v,u\in\tilde u$
there is a geodesic passing through $v_0,v,u$ and $d(v,u)=3\Delta+1$.
Then $\tilde G$ is a tree.
\end{lemma}

\proof
Let $\tilde d$ denote the distance in $\tilde G$, and
set $\tilde v_0=\{v_0\}\in\tilde V$, the equivalence class
of $v_0$.
It is clear that if $v\in\tilde v\in \tilde V$ then
$\tilde d(\tilde v,\tilde v_0)= (3\Delta+1)d(v,v_0)$.

Let $\tilde v\in \tilde V-\{\tilde v_0\}$.
By construction, $\tilde v$ has no neighbor $\tilde u$
with $\tilde d(\tilde v,\tilde v_0)=\tilde d(\tilde u,\tilde v_0)$.
We claim that there is exactly one neighbor $\tilde u$ of
$\tilde v$ such that
$\tilde d(\tilde v,\tilde v_0)=\tilde d(\tilde u,\tilde v_0)+1$.
Indeed, let $v_1,v_2\in\tilde v$, and suppose that
$d(v_1,v_2)\leq 2 \Delta$.
Let $s_1$ be a shortest path from $v_1$ to $v_0$, and
let $s_2$ be a shortest path from $v_2$ to $v_0$.
Let $u_1$ be the vertex on $s_1$ whose distance from $v_1$
is equal to $3\Delta+1$, and let $u_2$ be the vertex on
$s_2$ whose distance from $v_2$ is $3\Delta+1$.
Let $s_3$ be a shortest segment joining $v_1$ and $v_2$.
Because $G$ is $\Delta$-hyperbolic, there is a vertex $w$
in $s_3\cup s_2$ whose distance to $u_1$ is at most $\Delta$.
But the length of $s_3$ is at most $2\Delta$, and hence the
distance from $s_3$ to $u_1$ is at least $\Delta+1$.
Therefore, $w\notin s_3$ and $w\in s_2$.
Because $|d(w,v_0)-d(u_2,v_0)|=|d(w,v_0)-d(u_1,v_0)|\leq\Delta$,
and $w\in s_2$, we get $d(w,u_2)\leq \Delta$.
So
$d(u_1,u_2)\leq d(u_1,w)+d(w,u_2)\leq 2 \Delta$,
which implies that $u_1$ and $u_2$ belong to the same equivalence
class in $\tilde V$.
It now follows that $\tilde v$ has exactly one neighbor $\tilde u\in\tilde V$
with $\tilde d (\tilde v,\tilde v_0)=\tilde d (\tilde u,\tilde v_0)+1$.

Suppose that $\tilde G$ is not a tree.
Then there is a simple closed loop $\tilde C$ in $\tilde G$.
Consider a vertex $\tilde v$ in $\tilde C$, where
$\tilde d(\tilde v,\tilde v_0)$ is maximal.  The two
edges adjacent to $\tilde v$ in $\tilde C$ must connect
$\tilde v$ with neighbors $\tilde u$ with
$\tilde d (\tilde v,\tilde v_0)\geq \tilde d(\tilde u,\tilde v_0)$,
but $\tilde v$ has at most one such neighbor.  This contradicts
the existence of $\tilde C$, and hence $\tilde G$ is a tree.
\qed
\medskip

\noindent
{\bf Proof of Theorem~\ref{hypqit}}
Let $v_0$ be some vertex in $G$,
and let $V_1$, $\tilde G$ be as in Lemma~\ref{tree}.
For each $v\in V$ choose some geodesic $\sigma(v)$ from $v$ to $v_0$,
and let $v'$ be the vertex closest to $v$ in $\sigma(v)\cap V_1$.
Let $\pi v$ denote the equivalence
class of $v'$ in $\tilde V=V/\sim$.  It is an easy
exercise to check that the $\Delta$-hyperbolicity of $G$
implies that $\pi$ is a contraction.
On the other hand, it follows from Lemma~\ref{lemdiam}
that $\tilde d(\pi v,\pi u) \geq C^{-1} d(v,u) -C$
holds for some $C>0$ and every $v,u\in V$.
We conclude that $\pi$ is a rough-isometry.
Now lemma~\ref{tree} shows that $\tilde G$ is a tree,
and completes the proof of the theorem.
\qed

\begin{question}
Show that for any planar hyperbolic graph $G$,
$$
sep_G(n) < C\log n,
$$
for some $C(G)=C > 0$.
\end{question}

Note that this cannot be true in general if we omit planarity, by Proposition~\ref{hypspace}.

\section{More on regular maps}\label{s.reg}
Let $T$ denote the binary tree. The next result is shown to be true in \cite{BS}, even with a quasi-isometric embedding.

\begin{theorem}
\label{hrtt}
$\H^2 \_reg T \times T$.
\end{theorem}

\medskip

Using Theorem~\ref{TxT}
and the remark after its proof, the previous theorem implies that $sep_{\H^2 \times \H^2}(n) \asymp n/ \log n$,
because $T\times T\_reg \H^2 \times \H^2 \_reg  T \times T \times T \times T$.

\medskip
\begin{question}
\label{HTZ}
Is $\H^2 \_reg T \times \Z$?
\end{question}


There may be something special about the $n/ \log(n)$ separation function. The following problems are weaker than Question~\ref{linear}.

\begin{question}

\begin{enumerate}

\item
Suppose $G,H$ have separation functions $ \leq n/ \log n$.  Is it true
that $G \times H$ has separation function $\leq n/ \log n$?



\item Is it true that balls in infinite transitive graphs have
$$
\sep (B(r)) \leq
|B(r)|/ \log |B(r)|?
$$
\end{enumerate}
\end{question}
\medskip

The last question, if true, would answer negatively an old question of the first author, namely, if there is any Cayley graph where the balls form an expander family (see \cite{BK}). 
Without assuming transitivity, the answer to this question is clearly no: take some infinite tree with expanders put on each set of points at the same distance from the origin (as in \cite{BK}).


We have found that $sep_{\H^2 \times \R}(n)$ is at least $\sqrt{n \log n}$ (Corollary~\ref{calcul}).
This implies that
$\H^2 \times \R$ does not have a regular map to $\H^3$ which has separation only of
order $\sqrt n$ (Proposition~\ref{hypspace}).

\begin{question}
Is $\H^3 \_reg \H^2 \times \R$?
\end{question}
We doubt that the answer would be positive, but as we  said in Section~\ref{s.def}, we know only four ways to rule out regular maps:
separation, Dirichlet harmonic functions, asymptotic dimension and growth.  Neither does the job here.

\medskip



\section{Semi-regular maps} \label{s-reg}

In a private communication, Yehuda Shalom asked:
Is it true that a Cayley graph of an amenable group which is not
virtually cyclic admits a bilipschitz embedding of $\Z^2$ ?
Yehuda suggested that the lamplighter group $(\Z/2\Z)_2\wr \Z$ is a counterexample (see e.g. \cite{N} for a detailed definition).

Here is a proof not using any of the
arguments applied in the sections before but
rather using the notion of semi-regular maps defined
below. This will allow us to prove
that there is no regular map $\Z^2  \_reg  \mbox{Lamplighter}(\Z)$.
\medskip

\noindent
{\bf Definition.}
Call a map $f : X \rightarrow Y$   {\bf semi-regular}, denoted $X \rightarrow_{\bf s-reg} Y$,
if $f$ is Lipschitz and for
every $r$ there is a $c(r)< \infty$ such that for every
$y \in Y$ every connected component of $f^{-1}(B(y,r))$
has diameter at most $c(r)$.
\medskip

It is an easy exercise to verify that a composition of
semi-regular maps is semi-regular.
\medskip

Also, a regular map is semi-regular (we assume bounded
degree.  There are analogous definitions appropriate
for general metric spaces).
\medskip

Note that the canonical map $\mbox{Lamplighter}(\Z) \rightarrow \Z$
(location of the lamplighter) is semi-regular,
with $c(r)$ order $r2^r$.
\medskip

\begin{proposition}
There is no semi-regular map $\Z^2  \rightarrow_{\bf s-reg} \Z$.
\end{proposition}

By combining with the above, it follows that there is
no semi-regular map $\Z^2 \rightarrow_{\bf s-reg} \mbox{Lamplighter}(\Z)$.  In particular,
there is no regular map $\Z^2  \_reg  \mbox{Lamplighter}(\Z)$.
\medskip

\noindent
{\bf Proof:}
Let $g: Z^2 \rightarrow \Z$ be Lipschitz.
We  look at ``quasi-level sets" of $g$, and show that there
is some quasi-connected component of a quasi-level set that
is large (at this point, this is a vague statement;
the precise formulation will become clear shortly).
Take $k$ large, and consider the map $a(z) := \lfloor(z/k)\rfloor$,
$a: \Z-> \Z$.  Fix $k$ large enough so that $|a(g(z))-a(g(z'))| < 2$
if $z$ and $z'$ are within distance $2$ in $\Z^2$.  Let $X$ be the graph
obtained from $\Z^2$ by adding diagonals.  Given $r \in \Z$, let
$S(r)$ denote the set of connected components in $X$ of
$\{x \in X: a(g(z))= r\}$.  Let $S* := \cup_{r \in \Z} S(r)$.
This is a partition of $X$.
Note that if $S \in S(r)$ is finite, then there is a unique
$S' \in S(r+1)\cup S(r-1)$ that is adjacent to $S$ and ``surrounds"
it.  Let $r \in \Z$ be such that $S(r)$ is nonempty.  Let $S_0 \in S(r)$.
Inductively, define $S_{n+1}$ to be $S_n$ if $S_n$ is infinite,
and if not, let $S_{n+1}$ be some $S \in S*$ which surrounds $S_n$.
Clearly, $|S_n| \ge n$.  Hence $g$ is not semi-regular.
\qed
\medskip

Note: this also shows that there is no regular map
\newline
$\Z^2  \_reg  \mbox{Lamplighter}(\mbox{Lamplighter}(\Z))$.
\medskip

There is a more general way to see the proposition, which we will sketch in a later remark.

Before question~\ref{HTZ}, maybe consider,

\begin{question}
Does
$ \H^2 \rightarrow_{\bf s-reg} T \times \Z$?
\end{question}

\medskip

Gromov~\cite{G3} defined the notion of asymptotic dimension of a metric space,
which is defined as follows.
Say that $X$ has asymptotic dimension at most $D$ if
for every $s>0$ there is a partition (i.e., coloring)
$X = Y_1 \cup Y_2 \cup ... \cup Y_{D+1}$
such that each $Y_j$ can be partitioned into a collection of
pieces of bounded diameter (depending only on $s$) such that the distance between any
two of the pieces is at least $s$. See \cite{BD} for a nice survey on the asymptotic dimension.

Now, if $X$ maps semi-regularly into $\Z^d$, then the asymptotic
dimension of $X$ is bounded from above by the asymptotic dimension of $\Z^d$ (which is $d$, \cite{BD}). To see this, just pull back the colorings in the definition of asymptotic dimension using the semi-regular map.
By a similar argument one can show that there is no (semi-)regular map from $T\times T\times T$ to $T\times T$ (because the asymptotic dimension of the former is 3, and that of the latter is 2), and 
that there is no (semi-)regular map from $\R^d$ ($d\geq 3$) to $T\times T$ (even though the separation of $T\times T$ is larger). 

The question is whether this observation captures the partial
order among ``good" classes of spaces.  
To better understand the partial order given by semi-regular
embeddings, the following is a relevant question.

\begin{question}
Suppose that $G$ is a transitive (homogeneous) locally finite graph
with asymptotic dimension $d$.  Does it follow that $\Z^d \rightarrow_{\bf s-reg} G$
and $G \rightarrow_{\bf s-reg} \Z^d$ ?
\end{question}

\medskip
\noindent
{\bf Remark:}
We have just mentioned that the asymptotic dimension is monotone increasing under semi-regular maps. Hence there is no semi-regular map from $\Z^n$ to $\Z^m$ when $n>m$, giving an alternative proof for the above proposition.
\medskip

\noindent
{\bf Remark:}
Later Oded showed that  $\H^2$ and $\R^2$ are semi-regularly equivalent. The proof is unfortunately lost.

\medskip

{\bigskip\noindent\bf Acknowledgment}. Thanks to Gil Kalai, Mike Freedman and Yehuda Shalom for useful
discussions. Thanks to David Wilson for help with locating old emails and notes by Oded. Thanks to Sasha Sodin
for his useful comments on a previous version. We are grateful to Tatiana Smirnova-Nagnibeda and an anonymous referee for valuable suggestions.


\end{document}